\documentclass[10pt]{amsart}
\usepackage[cp1251]{inputenc}
\usepackage[english,russian]{babel}
\usepackage{amsmath,amsfonts,amssymb}
\usepackage[left=2cm]{geometry}
\newtheorem{lemma}{Lemma}
\newtheorem{theorem}{Theorem}

\newtheorem{corollary}{Corollary}
\newtheorem{prop}{Proposition}
\newcommand {\bproof }{{\par\medskip\noindent \bf Proof. }}
\newcommand {\eproof }{\hfill $\blacktriangle$ \\ \medskip}
\addto\captionsrussian{}

\begin{document}

\title[]{On completely regular codes with minimum eigenvalue in geometric graphs}
\author{I.Yu.Mogilnykh, K. V. Vorob'ev}

\thanks{\rm  This work was funded by the Russian Science
Foundation under grant 22-21-00135.}

\email{ivmog84@gmail.com, konstantin.vorobev@gmail.com}

\maketitle
\begin{quote}
{\begin{center} ABSTRACT \end{center}}

We prove that any completely regular code with minimum eigenvalue in any geometric graph $\Gamma$ corresponds to a
completely regular code in the clique graph of $\Gamma$. Studying the interrelation of these codes,  a complete characterization of the completely regular codes
in the Johnson graphs $J(n,w)$ with covering radius $w-1$ and strength $1$ is obtained.
In particular this result finishes a characterization of the completely regular codes in the Johnson graphs $J(n,3)$.
We also classify the completely regular codes of strength $1$ in the Johnson graphs $J(n,4)$ with only one case for the eigenvalues left open.
 \end{quote}
\medskip


\section{Introduction}

The geometric graphs are the distance-regular graphs, whose edges are parted into maximum cliques attaining Delsarte upper bound
 \cite{Godsil}. The Delsarte cliques are known to be completely regular codes in distance-regular graph.
Examples of geometric graphs include such classical distance-regular graphs as: Hamming, Johnson, Grassmann, dual-polar
graphs, etc. \cite{BCN}.



In recent papers \cite{MV}, \cite{GGV}, \cite{FI}, using various techniques, the classification problem for completely regular codes
in Hamming, Johnson and Grassmann graphs is reduced to the classification of codes or combinatorial structures in simplier graphs. In a sense, this paper continues this trend as we reformulate the existence problem for the completely regular codes with minimum eigenvalue in any geometric graph as that for completely regular codes in its clique graph.

For a survey on completely regular codes in Johnson graphs we refer to \cite{BRZ}.
In Johnson graphs, completely regular codes are tightly interrelated with important
classes of combinatorial $t$-designs such as Steiner triple and quadruple systems \cite{Mart},
\cite{AM}, \cite{Mog}. We refer to maximum $t$ such that a collection $D$ of fixed size blocks is a $t$-design as the {\it strength} of $D$.
The completely regular designs in $J(n,w)$ were classified in the following cases, see \cite{Mey}, \cite{M1}, \cite{GGV}, \cite{Vor}:

$ \bullet   $ strength $0$

$ \bullet   $ strength $1$ and minimum distance at least $2$

$ \bullet   $ strength $1$ and covering radius $1$.

Definitions and a review of the previous results are given in Section 2. A general approach is developed in Section 3.
 We show that any completely regular code $C$ with minimum eigenvalue in a geometric graph $\Gamma$  implies that
the set of maximum cliques having at least one vertex from $C$ is a completely regular code in the clique graph of $\Gamma$.
The proof is based on the fact that any Delsarte clique is a completely regular code in a distance-regular graph and the minimum eigenvalue of
the distance-regular graph is not its eigenvalue (see \cite[Proposition 2]{KMP} and \cite[Proposition 4.4.6]{BCN}).

The approach is implemented in Section 4,
where we give a characterization of the completely regular codes in the Johnson graphs $J(n,w)$ with covering radius $w-1$ and strength $1$.
This, combined with the results of work \cite{Mey} of Meyerowitz,
finish the characterization of completely regular codes of covering radius $w-1$ and any strength. 
We show that these codes for $w=3$ are necessarily a part of the series from the work \cite{GP} of Praeger and Godsil (see also \cite{M1}),  and there are no such codes for $w\geq 4$.
 This result in combination with \cite{Mey}, \cite{GGV}, \cite{GG}, \cite{Deh} closes the classification problem of the completely regular codes in $J(n,3)$. The completely regular codes of strength $2$ coincide with simple $2-(n,3,\lambda)$ designs \cite{Mart}.
Since these important combinatorial objects include Steiner triple systems, their complete classification up to isomorphism is far from being solved.
However, we have their classification up to parameter $\lambda$, see work \cite{Deh} of Dehon.

The classification results for completely regular codes in the Johnson graphs are obtained following the general idea for the geometric graphs. The Johnson graph $J(n,w-1)$ can be viewed as the clique graph of the geometric graph $J(n,w)$.
From results of Section 3, any completely regular code in the Johnson graph $J(n,w)$ with covering radius $w-1$ and strength $1$  
naturally corresponds to a completely regular code in $J(n,w-1)$ with covering radius $w-2$ and strength $1$. The incidence relations
between these codes impose constraints on their structure. In particular, once the code in $J(n,w-1)$ is known, one can reconstruct the corresponding code in $J(n,w)$.

We also show that all completely regular codes in $J(n,4)$ of strength $1$ are either part of known series \cite{M1}, \cite{GP} or their eigenvalues are $\theta_0(n,4),\theta_2(n,4),\theta_3(n,4)$. 
In the latter case using several local and global arguments we were able to find significant restrictions on the intersection arrays of such codes.


\section{Definitions and basic theory}
\subsection{Geometric graphs and their clique graphs}
A number $\theta$ is called an {\it eigenvalue} of a graph $\Gamma$ if there is a nonzero vector $v$ indexed by the vertices of $\Gamma$ such that
$$Av=\theta v,$$
where $A$ is the adjacency matrix of $\Gamma$. We call such vector $v$ a $\theta$-{\it eigenvector} of $\Gamma$.

A regular graph of diameter $d$ is called {\it distance-regular} if there is an array of integers $\{\beta_0,\ldots,\beta_{d-1};\gamma_1,\ldots,\gamma_{d}\}$, such that for any vertices $x$ and $y$ at distance $i$, $i\in\{0,\ldots,d\}$ there are
exactly $\beta_i$ neighbors of $y$ at distance $i+1$ from $x$ and $\gamma_i$  neighbors of $y$ at distance $i-1$  from $x$. The array of integers $\{\beta_0,\ldots,\beta_{d-1};\gamma_1,\ldots,\gamma_{d}\}$
is called the {\it intersection array} of the distance-regular graph $\Gamma$.
It is well known that any distance-regular graph of diameter $d$ has exactly $d+1$ distinct eigenvalues. We index these eigenvalues in descending order throughout the paper:
$$\theta_0=k> . . . > \theta_d,$$ where $k$ is the valency of $\Gamma$.
 If $K$ is a clique in a distance-regular graph $\Gamma$ of valency $k$ and diameter $d$, then \cite[Proposition 4.4]{BCN}
$$|K|\leq 1-\frac{k}{\theta_d}.$$

The bound was established by Delsarte \cite{Del}. A clique that attains the above upper bound is called {\it a Delsarte clique} of $\Gamma$. Any Delsarte clique is a completely regular code in $\Gamma$ \cite[Proposition 4.4.6]{BCN}.
If there is a set ${\mathcal K}$ of Delsarte cliques in a distance-regular graph such
that any edge of the graph is in a unique clique of ${\mathcal K}$, the graph is called
{\it geometric} with respect to ${\mathcal K}$ \cite{Godsil}. Note that there can be many ways to choose
the set ${\mathcal K}$. A generalization of the concepts of Delsarte clique graph
and geometric graph to non-distance-regular case is known as completely
regular clique graphs \cite{Suz}.

Let $\Gamma$ be a geometric graph with respect to a set of Delsarte cliques ${\mathcal K}$.
{\it The clique graph }of $\Gamma$ with respect to ${\mathcal K}$ is the graph $\Gamma'$
such that

$$ V(\Gamma')={\mathcal K}$$
and
$$ E(\Gamma')=\{(K,K'): K, K'\in {\mathcal K}, |K\cap K'|=1\}.$$

 The vertices of the {\it Johnson graph} $J(n,w)$ are all $w$-subsets of $\{1,\ldots,n\}$, where two $w$-subsets are adjacent if they meet in a $(w-1)$-subset. Since the Johnson graphs $J(n,w)$ and $J(n,n-w)$ are
isomorphic, throughout the text we assume that $w\leq n/2$. The Johnson graph $J(n,w)$, $w\leq n/2$ is distance-regular of diameter $w$ and its eigenvalues are
$$\theta_i(n,w)=(n-w-i)(w-i)-i,$$
for all $i\in \{0,\ldots,w\}$.

{\bf Example 1.}
Given a $(w-1)$-subset  $y$ of $\{1,\ldots,n\}$, one can consider the collection of all $w$-subsets of $\{1,\ldots,n\}$ containing $y$, which is a Delsarte clique in the Johnson graph $J(n,w)$ \cite{BCN}. Moreover, the graph $J(n,w)$ is geometric with respect to the set ${\mathcal K}$
of these cliques for all $(w-1)$-subsets $y$ of $\{1,\ldots,n\}$. The clique graph of $J(n,w)$ can be represented as the graph with $(w-1)$-subsets of $\{1,\ldots,n\}$ as vertices with edges being subsets meeting in a $(w-2)$-subset, i.e. it is
actually the Johnson graph $J(n,w-1)$. A vertex $x$ of $J(n,w)$ is in the clique $y$ of ${\mathcal K}$ if and only if the vertex $x$ (as $w$-subset) contains the vertex y (as $(w-1)$-subset). So we see that the incidence of the
vertices of the geometric graph $J(n,w)$ and its clique graph $J(n,w-1)$ is the ordinary subset inclusion relation. We also note that in a similar manner the Grassmann graph $J_q(n,w)$ can be viewed as a geometric graph whose clique graph is $J_q(n,w-1)$.

For a vector $u$ indexed by the vertices of the graph $\Gamma$, we define the vector $I(u)$ indexed by the cliques from a set ${\mathcal K}$ such that for $K\in {\mathcal K}$ we have $(I(u))_{K}=\sum\limits_{x\in K} u_x$.
This mapping a particular case of the following function. For a bipartite graph $\overline{\Gamma}$ and a vector u indexed by the vertices of one part, define $(I(u))_y=\sum\limits_{(x,y)\in E(\overline{\Gamma})}u_x$ for all vertices y from another part. 

\begin{theorem}\label{TAM}\cite[Theorem 1]{AM} Let $\overline{\Gamma}$ be a biregular bipartite graph with valencies $c$ and $c'$ and the halved
graphs $\Gamma$ and $\Gamma'$, $u$ be a $\theta$-eigenvector of $\Gamma$.
Let any pair of vertices of $V(\Gamma)$ at distance $2$ in $\overline{\Gamma}$ have exactly $m$ common neighbors and any pair of vertices of $V(\Gamma')$ at distance $2$ in $\overline{\Gamma}$ have exactly  $m'$  common neighbors.
If $u$ is a $\theta$-eigenvector of $\Gamma$ then the vector $I(u)$
 is a $\frac{c-c'+m\theta}{m'}$-eigenvector of $\Gamma'$ if $\theta\neq-\frac{c}{m}$ and all-zero vector otherwise.

\end{theorem}

\begin{corollary}\label{T1} Let $\Gamma$ be a geometric graph with respect to a set of Delsarte
cliques ${\mathcal K}$, $\Gamma'$ be the clique graph of $\Gamma$ with respect to ${\mathcal K}$. Let $k$ be the valency
of $\Gamma$, $\theta_d$ be its minimum eigenvalue and $u$ be a $\theta$-eigenvector of $\Gamma$. Then
the vector $I(u)$ is a $(\theta-\theta_d- 1 + \frac{k}{\theta_d})$-eigenvector of $\Gamma'$
if $\theta\neq \theta_d$ and $I(u)$ is
the all-zero vector otherwise.\end{corollary}
\bproof
Consider the bipartite graph whose parts are the vertices of $\Gamma$ and the cliques from ${\mathcal K}$ with incident vertices and cliques  being the edges.
 From the definition of geometric graph we see that such incidence graph is biregular with
valencies $-\theta_d, 1-\frac{k}{\theta_d}$ and any pair of its vertices at distance $2$ are adjacent
to a unique vertex. By Theorem \ref{TAM}  when $\theta\neq \theta_d$, for any $\theta$-eigenvector $u$ of one of the halved graphs $\Gamma$ of the biparitite graph, the vector $I(u)$ is always a  $(\theta-\theta_d- 1 + \frac{k}{\theta_d})$-eigenvector of the second halved graph $\Gamma'$.
Again by Theorem \ref{TAM} when $u$ is any $\theta_d$-eigenvector of $\Gamma$, the vector $I(u)$ is always the all-zero vector.

\eproof

\subsection{Completely regular codes}
Let C be {\it a code}, i.e. a collection of vertices in a regular graph $\Gamma$ and $\chi_C$ denote the characteristic vector of $C$ in $\Gamma$. A vertex $x$ is in $C_i$ if the minimum of
the distances between $x$ and the vertices of $C$ is $i$. The maximum of these
distances is called the {\it covering radius} of $C$ and is denoted by $\rho$. The {\it distance
partition} of the vertices of $C$ with respect to $C_0 = C$ is $\{C_i : i \in\{0,\ldots,\rho\}\}$.
A code $C$ is called {\it completely regular} if there are numbers $\alpha_0, \ldots,\alpha_\rho$, $\beta_0,\ldots,\beta_{\rho-1}$, $\gamma_1,\ldots,\gamma_{\rho}$ such that any vertex of $C_i$
is adjacent to exactly $\alpha_i$, $\beta_i$, $\gamma_i$ vertices of $C_{i}, C_{i+1}$ and $C_{i-1}$ respectively. Note that $\alpha_0,\ldots, \alpha_\rho$ can be obtained from the remaining constants and the valency of the graph.
The array $\{\beta_0,\ldots,\beta_{\rho-1};\gamma_1,\ldots, \gamma_{\rho}\}$  is called {\it the intersection array} of the completely
regular code $C$.

For a completely regular code C consider the following tridiagonal $(\rho + 1) \times (\rho + 1)$
matrix
$$\left(%
\begin{array}{cccccc}
\alpha_0 & \beta_0& 0&0 &\ldots &0\\
\gamma_1&\alpha_1 & \beta_1&0 &\ldots &0\\
 .& .&.&.&.&.\\
.& .&.&.&.&.\\
0&.&0& \gamma_{\rho-1}&\alpha_{\rho-1} & \beta_{\rho-1}\\

0&.&0&0& \gamma_{\rho}&\alpha_{\rho}\end{array}%
\right),$$
 which we call the {\it intersection matrix}
of the completely regular code $C$. The eigenvalues of this matrix are called
{\it the eigenvalues } of the completely regular code $C$. The set of the eigenvalues of a completely regular code $C$ is denoted by $Spec(C)$.
Due to the tridiagonal structure of the above matrix, it has exactly
$\rho + 1$ distinct eigenvalues and $\rho + 1$ different eigenvectors up to the multiplication by a scalar. {\it The opposite code} of a code $C$ of covering radius $\rho$ is the code $C_{\rho}$.
It is easy to see that $C$ is a completely regular code if and only if the opposite code $C_{\rho}$ is completely regular. Moreover, if $C$ is completely regular, the eigevalues of $C$ and $C_{\rho}$ coincide.

{\bf Remark 1.} Again due to tridiagonal structure, we note that the first element of any eigenvector of the intersection matrix is never zero and no two consecutive elements in any eigenvector are zeros.

\begin{theorem}\label{T2} (Lloyd's theorem)\cite[Theorem 4.5]{CDZ} Let $C$ be a completely regular code with covering radius $\rho$ in a regular graph $\Gamma$. Then $C$ has exactly $\rho+1$  eigenvalues, they are all simple and  they are eigenvalues of $\Gamma$.
\end{theorem}

\begin{prop} \label{P1} \cite[Lemma 1]{M} Let $C$ be a completely regular code in a regular
graph $\Gamma$. Then there is a unique decomposition $\chi_C=\sum\limits_{\theta\in Spec(C)}u_{\theta}$, where $u_{\theta}$ is a $\theta$-eigenvector of $\Gamma$.


\end{prop}

The following fact is well known and could be found, for example, in \cite[Proposition 2]{KMP}.
\begin{prop}\label{P2}  Let $\Gamma$ be a distance-regular graph, $\theta_d$ be its minimum eigenvalue, $K$ be a Delsarte clique of $\Gamma$.
Then $K$ is a completely regular code in $\Gamma$ with $Spec(K)=Spec(\Gamma)\setminus \theta_d$ and the vector $\chi_K$ is orthogonal to any $\theta_d$-eigenvector of $\Gamma$.
\end{prop}

We finish this subsection with the following result for completely regular codes in the Johnson graphs.

\begin{prop}\label{prop_strength}\cite[Corollary 2.1]{Mart}
Let the eigenvalues of a completely regular code $C$  with covering radius $\rho$ in the Johnson graph $J(n,w)$  be

$$\theta_0(n,w)>\theta_{i_1}(n,w)>\ldots>\theta_{i_{\rho}}(n,w).$$
Then
the strength of $C$ is $i_1-1$ and in particular $\rho\leq w-i_1+1$.
\end{prop}

\section{Completely regular codes in geometric graphs with minimum eigenvalue}

\begin{theorem}\label{T3} Let $\Gamma$ be a geometric graph with respect to a set of Delsarte
cliques ${\mathcal K}$ and the minimum eigenvalue $\theta_d$, $C$ be a completely regular code with covering radius $\rho$ such that
$\theta_d\in Spec(C)$. Then there are nonzero numbers $a_i$, $i\in\{0,\ldots,\rho-1\}$ that depend only on $i$ and the intersection arrays of $C$ and $\Gamma$ such that for any $K\in {\mathcal K}$, $K\subset C_i\cup C_{i+1}$, we have
that $|K\cap C_i|=a_i$.
\end{theorem}
\bproof
Let $v$ be a $\theta_d$-eigenvector of the intersection matrix of $C$. Let
$v^{\Gamma}$ be the vector indexed by the vertices of $\Gamma$ such that $v^{\Gamma}_x=v_i$ if $x$ is in $C_i$. It is well known that
$v^{\Gamma}$ is a $\theta_d$-eigenvector of $\Gamma$, see \cite[Theorem 4.5]{CDZ}, which, in particular, implies the Lloyd's theorem (Theorem \ref{T2}).

By definition of a geometric graph, for any $i\in\{0,\ldots, \rho-1\}$, any two
adjacent vertices of $C_i$ and $C_{i+1}$ are in a unique clique $K$ of ${\mathcal K}$. Note that
since $K$ is a clique and $C_0,\ldots,C_{\rho}$ is a distance partition, $K$ is a
subset of $C_i \cup C_{i+1}$.
By Proposition \ref{P2}, the $\theta_d$-eigenvector $v^{\Gamma}$ is orthogonal to $\chi_K$. So, for all $i\in\{0,\ldots, \rho-1\}$ there is $K\in {\mathcal K}$:
\begin{equation}
\label{eq1}
v_i|K\cap C_i|+v_{i+1}(|K|-|K\cap C_i|) = 0. \end{equation}

By the choice of cliques $K$'s any such $K$ contains an edge between $C_i$ and
$C_{i+1}$ and we have that $|K\cap C_i|\neq 0, |K|$. In view of Remark 1, we see that $v_i$ and
$v_{i+1}$ cannot be both zeros, so equality (\ref{eq1}) implies that $v_i\neq 0$ for all
$i$ from $\{0,\ldots,\rho\}$. Moreover, $v_i\neq v_{i+1}$ because otherwise from (\ref{eq1}) $K$ is empty. 
Suppose that there is $K\in {\mathcal K}$  such that $K \subset C_i$ for some $i$. Then
again by orthogonality of $v^{\Gamma}$ and $\chi_K$, we have that $|K|v_i = 0$, which contradicts $v_i\neq
0$. We conclude that any clique $K\in {\mathcal K }$ fulfills (\ref{eq1}) with $|K\cap C_i|\neq 0,|K|$ for some $i$.

We note that the eigenvector $v$ of the intersection matrix of $C$ is unique up to the multiplication by a scalar.
 This combined with fact that  the size of a Delsarte clique $K$ is uniquely determined by the intersection array of $\Gamma$ and (\ref{eq1})
imply that  the size of $K\cap C_i$ is uniquely determined by the intersection arrays of $C$ and $\Gamma$.
 \eproof

For a completely regular code $C$ in a geometric graph $\Gamma$ with respect to a set ${\mathcal K}$ of Delsarte cliques  consider the following collection of cliques:

$$C' =\{K\in {\mathcal K}: K\cap C \neq \emptyset\}.$$

In the next lemmas we work towards showing that $C'$ is a completely regular code in the clique graph of the geometric graph $\Gamma$.

\begin{lemma}\label{L1} Let $\Gamma$ be a geometric graph with respect to a set ${\mathcal K}$ of Delsarte
cliques and $\theta_d$ be the minimum eigenvalue of $\Gamma$ and let $\Gamma'$ be the clique graph of $\Gamma$ with respect to $\mathcal K$. If $C$ is a completely
regular code with covering radius $\rho$ in $\Gamma$ and eigenvalue $\theta_d$, then the distance partition with respect to the code $C'$ in $\Gamma'$ is $C'_i=\{K\in {\mathcal K}:K\subset C_i\cup C_{i+1}\}$, $i\in \{0,\ldots,\rho-1\}$.
\end{lemma}
\bproof The proof of Theorem \ref{T3} implies that for any $i\in\{0,\ldots, \rho-1\}$, the set $\{K\in{\mathcal K}:K\subset C_i\cup C_{i+1}\}$ actually consists of the cliques
in ${\mathcal K}$ with edges between $C_i$ and $C_{i+1}$.
Since $C_0,\ldots, C_{\rho}$ is the distance partition with respect to the completely regular code $C$, any edge between $C_i$ and $C_{i+1}$
is incident to an edge between $C_{i-1}$ and $C_{i}$, for any $i$ such that $1\leq i\leq \rho-1$.
In view of the above reference to the proof of Theorem \ref{T3}, any clique in
$\{K\in {\mathcal K}: K\subset C_i\cup C_{i+1}\}$ meets a clique from $\{K\in{\mathcal K}: K\subset C_{i-1}\cup C_i\}$  in
a vertex. Equivalently, in terms of adjacency in the clique graph $\Gamma'$ a vertex of
$\{K\in {\mathcal K}: K\subset C_i \cup C_{i+1}\}$ has a neighbor in $\{K\in {\mathcal K}: K\subset C_{i-1}\cup C_i\}$. We obtain the required.
\eproof

\begin{lemma}\label{L2} Let the conditions of Lemma \ref{L1} hold and $x$ be a vertex of $C_i$ for
 $i\in \{0,\ldots,\rho\}$. Then $x$ is contained only in cliques from $C'_{i-1}$ and $C'_i$ and the
number of cliques from $C'_{i-1}$ containing $x$ is $b_i$, where $b_i$ depends only on $i$ and the intersection
arrays of $C$ and $\Gamma$. Moreover, we have that:

$$b_0=0, b_{\rho}=-\theta_d\mbox{ and }0<b_i<-\theta_d,$$ for $i\in \{1,\ldots,\rho-1\}$,
where $\theta_d$ is the minimum eigenvalue of $\Gamma$.
\end{lemma}
\bproof
By the definition of the completely regular code $C$, since $x$ is in $C_i$, $x$
has $\gamma_i$ neighbors in $C_{i-1}$, where
 $\gamma_i$ is an element of the intersection array
of $C$.

 On the other hand, by the definition of a geometric graph, the vertex $x$ and any of its $\gamma_{i}$ neighbors in $C_{i-1}$ are contained in only one clique from ${\mathcal K}$.
We note that all these cliques are in $C'_{i-1}=\{K\in {\mathcal K}:K\subset C_{i-1}\cup C_i\}$ and let $b_i$ be their number. By Theorem \ref{T3} if a clique belongs to $C'_{i-1}$, it has $a_{i-1}$ vertices of $C_{i-1}$.
We see that $x$ is in $b_i=\frac{\gamma_i}{a_{i-1}}$ cliques from $C'_{i-1}$. By the proof of  Theorem \ref{T3}, any another clique containing a vertex $x$ of $C_i$, contains vertices of both $C_i$ and $C_{i+1}$, i.e. this clique is in $C'_i$.

Since there is no $C'_{-1}$ and $C'_{\rho}$, we see that $b_0=0$ and $b_{\rho}$ equals the number of the cliques in ${\mathcal K}$ containing a fixed vertex. The number $b_{\rho}$ is $\frac{k}{|K|-1}=-\theta_d$, where $|K|$ is the size of a Delsarte clique $K$ in $\Gamma$ and $k$ is the valency of the graph $\Gamma$.
We note that for $i \in \{1,\ldots,\rho-1\}$, $b_i$ cannot be $0$ or $-\theta_d$. Indeed, any vertex in $C_i$ has neighbors
in $C_{i-1}$ and $C_{i+1}$ because $C$ is a completely regular code. Since the edgeset of $\Gamma$ is parted into cliques of ${\mathcal K}$, $x$ belongs to at least one clique
 contained in  $C_{i-1}\cup C_{i}$ and at least one clique  contained in  $C_{i}\cup C_{i+1}$. In other words, taking into account the expression for $C'_{i}$'s, $x$ belongs to at least one clique in $C'_{i-1}$ and $C'_{i}$ and $0<b_i<-\theta_d$.

\eproof

\begin{theorem} Let $\Gamma$ be a geometric graph with respect to a set of Delsarte
cliques ${\mathcal K}$, valency $k$ and minimum eigenvalue $\theta_d$. Let $\Gamma'$ be the clique graph of $\Gamma$ with
respect to ${\mathcal K}$. If $C$ is a completely regular code in $\Gamma$, $\theta_d\in Spec(C)$ then
$\{K\in{\mathcal K}:K\cap C\neq \emptyset\}$ is a completely regular code in $\Gamma'$ and $Spec(C') =\{\theta-\theta_d-1+\frac{k}{\theta_d}:\theta\in Spec(C),\theta\neq \theta_d\}$.

\end{theorem}

\bproof
By Lemma \ref{L1}, for the code $C'=\{K\in{\mathcal K}:K\cap C\neq \emptyset\}$, $C'_i= \{K\in {\mathcal K}:K\subset C_i\cup C_{i+1}\}$, $i \in \{0, \ldots , \rho - 1\}$ is the distance partition of the vertices of the clique graph of $\Gamma$
 with
respect to the code $C'$. 

For any $i\in\{0,\ldots,\rho\}$ we consider a clique $K\in C'_i$ and find the number of cliques in $C'_{i-1}$, $C'_i$ and $C'_{i+1}$ meeting $K$.
By Theorem \ref{T3}, $|K \cap C_i |=a_i$ is the same for all $K \in C'_i$.
 By Lemma \ref{L2}, any vertex in $K \cap C_i$ is contained in $b_i$ cliques from $C'_{i-1}$ and
$-\theta_d -b_i-1$ cliques from
$C'_{i}$ different from $K$. Again by Lemma \ref{L2}, any vertex of $K \cap C_{i+1}$
is contained in $b_{i+1}-1$ cliques from $C'_i$ different from $K$ and $-\theta_d-b_{i+1}$ cliques from $C'_{i+1}$.
Since any two cliques  in ${\mathcal K}$ either meet in a unique vertex or are disjoint, the numbers
of the cliques from $C'_{i-1}$, $C'_i$ , $C'_{i+1}$ meeting $K$ are $a_ib_{i}$, $(-\theta_d-b_i-1)a_i+(b_{i+1}-1)(1-\frac{k}{\theta_d}-a_i)$ and $(-\theta_d-b_{i+1})(1-\frac{k}{\theta_d}-a_i)$ i.e. $C'$ is a completely regular code in $\Gamma'$.

In order to find the spectra of $C'$, we now yield a relationship between the characteristic vectors of $C$ and $C'$ in the vertex sets of $\Gamma$ and
$\Gamma'$ respectively.
For this purpose we consider the linear transform $I$, defined before Theorem \ref{TAM}. The transform $I$ of a vector $u$ indexed by the vertices of $\Gamma$ is defined as follows:
$$(I(u))_K=\sum\limits_{x\in K}u_x.$$
We take $u$ to be $\chi_C$. If the clique $K$ is not in $C'$, by definition of $C'$, $K$ does not contain vertices of $C$, so $(I(\chi_C))_K=\sum\limits_{x\in K} (\chi_C)_x=0$. If $K$ is in $C'$, then it contains exactly $a_0$ vertices of $C$ by Theorem \ref{T3},
so $(I(\chi_C))_K=a_0$. We see that

\begin{equation}
\label{eq2}
a_0\chi_{C'}=I(\chi_C).
\end{equation}

By Proposition \ref{P1} we have the decomposition
$$\chi_C=\sum\limits_{\theta\in Spec(C)} u_{\theta},$$ where $u_{\theta}$ is a $\theta$-eigenvector of $\Gamma$. In view of Corollary \ref{T1}, the vector $I(u_{\theta})$ is a $(\theta-\theta_d-1+\frac{k}{\theta_d})$-eigenvector of $\Gamma'$ if $\theta\neq \theta_d$ and all-zero vector if $\theta=\theta_d$. So, from (\ref{eq2}) we have that

$$a_0\chi_{C'}=\sum\limits_{\theta\in Spec(C),\theta\neq \theta_d}I(u_{\theta}).$$

We see that the characteristic vector $\chi_{C'}$ of the completely regular code $C'$ with covering radius $\rho-1$ is a sum of $\rho$ eigenvectors with pairwise distinct eigenvalues $\{\theta-\theta_d-1+\frac{k}{\theta_d}:\theta\in Spec(C),\theta\neq \theta_d\}$. This combined with
 Proposition \ref{P1} applied to the completely regular code $C'$ implies that $\{\theta-\theta_d-1+\frac{k}{\theta_d}:\theta\in Spec(C),\theta\neq \theta_d\}$ is the set of the eigenvalues of $C'$.
\eproof

\begin{corollary}\label{C1}
Let $C$ be a completely regular code in the Johnson graph $J(n,w)$
with covering radius $\rho$ and $\theta_w(n,w)\in Spec(C)$. Then the following holds for the code $C'=\{y:|y|=w-1,\exists x \in C\mbox{ such that } y\subset x\}$.

1. The code $C'$ is completely regular in $J(n,w-1)$ with covering radius $\rho-1$ and $Spec(C')=\{\theta_i(n,w-1):\theta_i(n,w)\in Spec(C), i\neq w\}$.
In particular the strengths of $C$ and $C'$ coincide.

2. There are constants $b_i$, $i=0,\ldots,\rho$ such that a vertex $x$ is in $C_i$ if and only if $\{y: |y|=w-1, y\subset x\}$ contains exactly $b_i$ vertices from $C'_{i-1}$ and $w-b_i$ vertices of $C'_i$. Moreover, we have $b_0=0$ and $b_{\rho}=w$ and
$0<b_i< w$ for $i \in \{1,\ldots,\rho-1\}$.
\end{corollary}
\bproof
We view the clique graph of the geometric graph $J(n,w)$ as the graph $J(n,w-1)$, where the vertices of $J(n,w-1)$ represent the Delsarte cliques
and vertex-clique inclusion relation between $J(n,w)$ and $J(n,w-1)$ is the inclusion relation of $(w-1)$-subset (clique) into $w$-subset (vertex). The complete regularity for $C'$ follows from Theorem \ref{T5}, whereas the existence of constants $b_i$'s is proved in Lemma \ref{L2}.
\eproof

\section{Characterization of  completely regular codes of strength $1$  in Johnson graphs}
\subsection{Constructions of completely regular codes of strength $1$}

In below we give the list of known compelely regular codes of strength $1$ in the Johnson graphs \cite{AM2}, \cite{GP}, \cite{M1}.
Let $Y_1, \ldots, Y_q$ be a partition of  $\{1,\ldots,n\}$
into $q$ subsets of same size $p$: $|Y_i|=p$ for $i=1,\ldots, q$.

GP.1 $p=2$ and odd $w$, $w\geq 5$. $C=\{x \mbox{ is a union of } \frac{w-1}{2} \mbox{ } Y_i \mbox{ 's  and one disjoint point}\}$

GP.1' $p=2$ and $w=3$. $C=\{x \mbox{ is a union of one } \mbox{ } Y_i \mbox{ for some }  i  \mbox{ and one disjoint point}\}$, $\rho=1$,
$Spec(C)=\{\theta_0(n,3),\theta_2(n,3)\}$.

GP.2 $p,q\geq 3$ and $w=3$. $C=\{x:\exists i, x\subseteq Y_i\}$. $\rho=2$,  $Spec(C)=\{\theta_0(n,3),\theta_2(n,3),\theta_3(n,3)\}$.

GP.3 $q=2$ and $w=3$. $C=\{x:\exists i, x\subseteq Y_i\}$; $\rho=1$, $Spec(C)=\{\theta_0(n,3),\theta_2(n,3)\}$.

GP.3' $q=2$ and $w=4$. $C=\{x:\exists i, x\subseteq Y_i\}$; $\rho=2 $, $Spec(C)=\{\theta_0(n,4),\theta_2(n,4),\theta_4(n,4)\}$.

GP.3'' $q=2$ and any $w$, $w\geq 5$. $C=\{x:\exists i, x\subseteq Y_i\}$; $\rho=\lfloor \frac{w}{2}\rfloor $.

GP.4 $p\geq 3$, $q\geq 2$ and $w=2$. $C=\{x:\exists i, x\subseteq Y_i\}$. $\rho=1$, $Spec(C)=\{\theta_0(n,2),\theta_2(n,2)\}$.

GP.5  $p=2$, $q\geq 2$ and $w=2$. $C=\{x:\exists i, x\subseteq Y_i\}$. $\rho=1$, $Spec(C)=\{\theta_0(n,2),\theta_2(n,2)\}$.

For the following two series we take $Y_i=\{2i-1,2i\}$, $i\in\{1,\ldots,\frac{n}{2}\}$ and even $n$.

M.1 Even $w$, $C=\{x:x \mbox{ is a union of }\frac{w}{2}\mbox{ } Y_i's\}$, $\rho=\frac{w}{2}$.

M.2 $w=3$, $C=\{x: |x\cap Y_i|\leq 1 \mbox{ for all } i \in \{1,\ldots,\frac{n}{2}\}, x \mbox{ contains both even and odd numbers}\}$,  $\rho=1$, $Spec(C)=\{\theta_0(n,3),\theta_2(n,3)\}$.

\subsection{Completely regular codes of strength $1$ in Johnson graphs}
We start with the following description of completely regular codes in Johnson graphs $J(n,2)$ with strength $1$.
\begin{prop}
If $C$ is a completely regular code of strength $1$ in $J(n,2)$, $C$ is a simple $1$-$(n,2,\lambda)$-design, where $\lambda$ or $n$ are even and vice versa.
\end{prop}
\bproof
The completely regular codes in $J(n,w)$ of strength $w-1$ coincide with the class of simple $(w-1)-(n,w,\lambda)$-designs \cite{Mart}.
It is well known that such objects exist if and only if $n$ or $\lambda$ are even.
\eproof

Proposition above and Theorem below complete the characterization of completely regular codes in $J(n,w)$ of strength $1$ and covering radius $w-1$.

\begin{theorem}
\label{T5}
Up to the opposite code, the following holds:

1. For any $n$ if $C$ is a completely regular code in $J(n,3)$ with covering radius $2$ and strength $1$ then the code $C'=\{x:|x|=2, \exists y\in C, x\subseteq y\}$ is the completely regular code (GP.4) in $J(n,2)$.

2. For any $n$, the only completely regular codes in $J(n,3)$ with covering radius $2$ and strength $1$  are codes obtained by construction (GP.2).

3. For any $w\geq 4$ there are no completely regular codes in $J(n,w)$ with covering radius $w-1$ and strength $1$.
\end{theorem}
\bproof
A completely regular code $C$ in $J(n,w)$ with covering radius $w-1$ by Theorem \ref{T2} has $w$ pairwise distinct eigenvalues, one of which is $\theta_0(n,w)$. If the strength of the code $C$ is 1, then by Proposition \ref{prop_strength} the eigenvalues of all completely
regular codes in the proof are $\theta_0(n,w)\cup \{\theta_i(n,w):i \in \{2,\ldots,w\}\}$.

1. Let $C$ be a completely regular code in $J(n,3)$ with covering radius $2$ and strength $1$. Then from the remark above the eigenvalues of $C$ are $\theta_0(n,3)$ , $\theta_2(n,3)$, $\theta_3(n,3)$. By Corollary \ref{C1} we see that
$C'=\{x:|x|=2, \exists y\in C, x\subseteq y\}$ is a completely regular code in $J(n,2)$ with $\rho=1$, strength $1$ and eigenvalue $\theta_2(n,2)$. By the second statement of Corollary \ref{C1}, there are constants $b_0$, $b_1$, $b_2$, such that
\begin{gather}
b_0=0<b_1< b_2=3\mbox{ and }\nonumber\\
\{j,l,m\}\in C_i\Leftrightarrow b_i \mbox{ vertices of }\{j,m\}, \{m,l\}, \{j,l\}
\mbox{ are in }C'_{i-1}\mbox{ and }3-b_i\mbox{ are in }C'_i\label{eqT5}
\end{gather}
Considering (\ref{eqT5}) when  $\{j,l,m\}$
is in $C$, since $b_0=0$
  we see that  pairwise adjacent vertices in $J(n,2)$ $\{j,m\}$, $\{m,l\}$, $\{j,l\}$  belong to $C'$. Consider
the number $\alpha_0$ of the vertices of the completely regular code $C'$ adjacent to a fixed vertex of $C'$.
Since $C'$ is not independent we see that $\alpha_0\neq 0$.

Taking the opposite code of $C'$ if necessary, we fix $b_1=1$ in what follows.
Now if $\{j,m\}$ and $\{m,l\}$ are in $C'$, the values for $b_i$'s above and (\ref{eqT5}) imply that the vertex $\{j,l\}$ is also in $C'$, because otherwise exactly two of vertices $\{j,m\}$, $\{m,l\}$, $\{j,l\}$
are in $C'$, i.e. $b_1$ is $2$, a contradiction.

For a vertex $\{j,m\}$ in $C'$, consider its neighbors from $C'$. The consideration above implies that if $\{m,l\}$ is in $C'$, the vertex $\{j,l\}$ is also in $C'$.
We see that the neighbors of $\{j,m\}$ belonging to $C'$ are vertices $\{j,l_s\}$, $\{m,l_s\}$, for $s=1,\ldots,\alpha_0/2$.

We show that for distinct $s$ and $t$, the vertex $\{l_s,l_t\}$ is in $C'$. Indeed $\{j,l_s\}$ and $\{j,l_t\}$ are in $C'$, forcing $\{l_s,l_t\}$ to be in $C'$.
We conclude that $C'$ is actually the union of all $2$-subsets of some $(\alpha_0/2+2)$-subsets.  Note that by the argument after (\ref{eqT5}),
 we have $\alpha_0\geq 2$. We now show that any two different such  $(\alpha_0/2+2)$-subsets are disjoint. Suppose the opposite. Let the $2$-subsets of $Y_1$ and $Y_2$ be in $C'$, $|Y_1|=|Y_2|=\alpha_0/2+2$ and $\alpha_0/2+2>|Y_1\cap Y_2|\geq 1$. Take an element $j\in Y_1\setminus Y_2$, and an element $l\in Y_1\cap Y_2$, which exists due to $|Y_1\cap Y_2|\geq 1$.
A vertex $\{j,l\}$ of $J(n,2)$ is in $C'$, because it is a subset of $Y_1$.
It is adjacent to $\alpha_0$ $2$-subsets of $Y_1$. Moreover it it adjacent to $\alpha_0/2+2-|Y_2\cap Y_1|$ $2$-subsets of $Y_2$ that are not 2-subsets of $Y_1$. These are obtained by replacing $j$ in $\{j,l\}$ with any element from $Y_2\setminus Y_1$. We see that $\{j,l\}$ is adjacent to at least
$3\alpha_0/2+2-|Y_2\cap Y_1|>\alpha_0$ vertices of $C'$ because $Y_1$ and $Y_2$ are distinct subsets of size $\alpha_0/2+2$ and $\alpha_0/2+2>|Y_1\cap Y_2|$.
This contradicts the fact that according to the definition of the intersection array of the completely regular code $C'$ the number of the vertices of $C'$ being adjacent to a fixed vertex of $C'$ is $\alpha_0$.

Thus the code $C'$ is the union of the $2$-subsets of disjoint $(\alpha_0/2+2)$-subsets. If these  $(\alpha_0/2+2)$-subsets do not partition $\{1,\ldots,n\}$, there is an element $j$ of $\{1,\ldots,n\}$ such that
$\{j,l\}$ is not a vertex of $C'$ for any $l$.
We see that $C$ is not a $1$-design, a contradiction.

From the considerations in the begining of the proof $\alpha_0$ is nonzero, so $|Y_i|\geq 3$.

2. We have shown that $C'$ is the union of the $2$-subsets of disjoint subsets $\{Y_1,\ldots,Y_q\}$ of a fixed size $p$, $p=\alpha_0/2+2\geq 3$.  For such $C'$ we reconstruct the initial code $C$ in $J(n,3)$. By the second statement of Corollary \ref{C1} and taking
into account that $b_0=0$, a $3$-subset $\{j,l,m\}$ is a vertex of $C$ if and only if all its $2$-subsets are in $C'$. From the obtained description of $C'$, we see that this is the case only when there is $i$: $\{j,l,m\}$ is a subset of $Y_i$, i.e. the code $C$ is actually (GP.2).

3. We now show that there are no completely regular codes with $\rho=3$ and strength $1$ in $J(n,4)$. Suppose the opposite and let $D$ be such code and $D'=\{x:|x|=3, \exists y \in D, x\subseteq y\}$. Then by Corollary \ref{C1} the code $D'$ is completely regular
in $J(n,3)$ with covering radius $2$ and strength $1$. By the shown in the second statement of current Theorem, $D'$ is  the union of the $3$-subsets of disjoint subsets  $\{Y_1,\ldots,Y_q\}$ of fixed size $p, p\geq 3$. Moreover, again by the second statement of Corollary \ref{C1},
a $4$-subset of $\{1,\ldots,n\}$ is a vertex of $D$ if and only if all its $3$-subsets are in $D'$. Therefore, the code $D$ consists of the $4$-subsets of all $Y_i$'s. Further considerations are for the code $D$ only. Now we have two options. If there are only two such $Y_i$'s, $D$ is (GP.3') and therefore has covering radius $2$,
which contradicts that the covering radius of $D$ is $3$.   If there are at least $3$ such $Y_i$'s we show that $D$ is not completely regular. Indeed, it is not hard to see that the distance partition with respect to $D$
is such that:

$$D_1=\{x:|x|=4, \exists i\mbox{ } |x\cap Y_i|=3\} $$
and
$$D_2=\{x:|x|=4, \exists i\mbox{ } |x\cap Y_i|=2\}. $$
Let the subsets $Y_1$, $Y_2$ and $Y_3$ be such that $\{1,2\}\subset Y_1$, $\{3,4\}\subset Y_2$, $5\in Y_3$.
The subsets $\{1,2,3,4\}$, $\{1,2,3,5\}$ are both in $D_2$. However, in view of the description of $D_1$ and $D_2$ above, the vertex $\{1,2,3,4\}\in D_2$ is adjacent to $4(|Y_1|-2)=4(p-2)$ vertices of $D_1$ in the graph $J(n,4)$, whereas $\{1,2,3,5\}\in D_2$ is adjacent to
$2(|Y_1|-2)=2(p-2)$ vertices of $D_1$, i.e. $D$ is not completely regular. We conclude that there are no completely regular codes in $J(n,4)$ with $\rho=3$ and strength $1$.

By induction hypothesis on $w$, $w\geq 5$, we have that there are no completely regular codes in $J(n,w-1)$, covering radius $w-2$ and strength $1$. If a completely regular code with covering radius $w-1$ and strength $1$ exists in $J(n,w)$, then by Corollary \ref{C1} there is a completely
regular code in $J(n,w-1)$ with $\rho=w-2$ and strength $1$, a contradiction.

\eproof


According to Proposition \ref{prop_strength}, any completely regular code of covering radius $w-1$ has strength $0$ or $1$. The codes of strength $0$ were completely characterized in \cite{Mey}. From
Theorem \ref{T5} we have the following. 
\begin{corollary}
\label{CoroJnall}
Let $C$ be a completely regular code in $J(n,w)$ of covering radius $w-1$. Then it either has strength 1 and is described in Theorem 5 or it has strength 0 and characterized in \cite{Mey}.

\end{corollary}

\begin{corollary}
\label{CoroJn3}
Let $C$ be a completely regular code in $J(n,3)$. Then up to opposite code the following holds.
\begin{enumerate}
    \item \cite{Mey} If $C$ has strength $0$ then $C$ is  $\{x:I \subseteq x\}$ for any $I\subset\{1,\ldots,n\}$, $|I|\leq 3 $
or $\{x:x \subset I\}$,  for any $I\subseteq\{1,\ldots,n\}$, $|I|\geq 4 $.
    \item \cite{GGV} If $C$ has strength $1$ and covering radius $1$ then $n$ is even. Moreover, for $n=6$ the code $C$ is the union of arbitrary collection of antipodal classes in $J(6,3)$;
for $n\geq 8$  and $n\neq 10$, the code $C$ is (GP.1') for $w=3$, (GP.3) or (M.2), for n=10 the code is (GP.1') for w=3 , (GP.3) or has the same intersection array as (M.2).
    \item  If the code $C$ has strength $1$ and covering radius $2$ then $C$ coincides with (GP.2).
    \item \cite{Mart}, \cite{Deh} If the code $C$ has strength $2$ then $C$ is a simple $2-(n,3,\lambda)$-design with $\lambda n(n-1) = 0 (\mbox{mod }6)$ and $\lambda (n-1) = 0 (\mbox{mod }2)$ and vice versa.
\end{enumerate}
\end{corollary}
\bproof
Any completely regular code in $J(n,3)$ has strength $t$, for some $t$, $0\leq t\leq 2$  and covering radius $\rho$. By Proposition \ref{prop_strength}
 we have $\rho\leq 3-t$.
All completely regular codes in $J(n,3)$ for $t=0$ were characterized in \cite{Mey}. The completely regular codes in $J(n,3)$ for $t=1$ and $\rho=1$ are codes (GP.1') for $w=3$, (GP.3) or (M.2)  for even $n\geq 12$ and any odd $n$ \cite{GGV}.
The completely regular codes in $J(6,3)$ with $t=1$ and $\rho=1$ were characterized in \cite{AM2} in terms of antipodal classes.
 By computer search we verified that any completely regular code with $t=1$ and $\rho=1$ in $J(8,3)$ is (GP.1') for $w=3$, (GP.3) or (M.2). Moreover, any completely regular code with $t=1$ and $\rho=1$ in $J(10,3)$
is (GP.1') for $w=3$, (GP.3) or has the same intersection array as (M.2) (there are several such nonisomorphic codes).

By Proposition \ref{prop_strength} the codes with $t=1$ and $\rho=2$ in $J(n,3)$ necessarily have eigenvalues $\theta_0(n,3)$, $\theta_2(n,3)$, $\theta_3(n,3)$ and therefore are described in Theorem \ref{T5}. The codes in $J(n,3)$ of strength $2$ coincide with the class of simple $2-(n,3,\lambda)$-designs \cite{Mart}.
The classification of these designs
up to parameter $\lambda$ was obtained in \cite{Deh}: a simple $2$-$(n,3,\lambda)$-designs exists if and only if $\lambda n(n-1) = 0 (\mbox{mod }6)$ and $\lambda (n-1) = 0 (\mbox{mod }2)$.
\eproof

Now we can apply our technique based on Delsarte cliques for completely regular codes in the Johnson graphs $J(n,4)$.

\begin{theorem}\label{T5.5}
	Let $C$ be a completely regular code of strength $1$ in $J(n,4)$, $\theta_4(n,4)\in Spec(C)$. Then $C$ has covering radius $2$, $Spec(C)=\{\theta_0(n,4),\theta_2(n,4), \theta_4(n,4) \}$ and $C$ is (GP.3') or (M.1) for $w=4$.

%

\end{theorem}

\bproof

According to Proposition \ref{prop_strength} we have the following alternatives for any completely regular code $C$ of strength $1$ in $J(n,4)$ with $\theta_4(n,4)\in Spec(C)$:

\begin{enumerate}
	\item  $C$ has covering radius $2$, $Spec(C)=\{\theta_0(n,4),\theta_2(n,4), \theta_4(n,4) \}$.
	\item  $C$ has covering radius $3$ and $Spec(C)=\{\theta_0(n,4),\theta_2(n,4),\theta_3(n,4), \theta_4(n,4) \}$. Such codes do not exist by Theorem \ref{T5}.
\end{enumerate}

We proceed by considering the first case in more details. Let $C'$ be $\{y:|y|=3,\exists x \in C, y\subset x\}$. In view of Corollary \ref{C1} the code $C'$ has covering radius $1$ and $Spec(C')=\{\theta_0(n,3), \theta_2(n,3)\}$.
By \cite{GGV}, see Corollary \ref{CoroJn3}, $C'$ is one of the codes (GP.1'), (GP.3), (M.2) or codes with the  same intersection matrix as (M.2) for n=10. Let $Y_i$ be $\{2i-1,2i\}$ for $i=1,\ldots, \frac{n}{2}$.
By Corollary \ref{C1} one can reconstruct the code $C$ from $C'$ using constants $b_i$'s such that
$$b_0=0<b_1<b_2=4,$$
where $x\in C_i$ if and only if $x$ is contains $b_{i}$ vertices of $C'_{i-1}$ and $4-b_{i}$ vertices of $C'_{i}$.

Let $C'$ be $\{Y_i\cup j:i\in \{1,\ldots, n/2\}, \mbox{ for all }j \notin Y_i\}$, i.e. $C'$ is (GP.1').
We see that the only $4$-subsets of $\{1,\ldots,n\}$ such that any its $3$-subsets is in $C'$  are the unions of $Y_i$'s, i.e.

$$C=\{Y_i\cup Y_j: i\neq j\}, $$

which is exactly the code (M.1) for $w=4$.

If $C'$ is (M.2), we prove that $C$ does not exist by showing that  $b_1$ is inconsistent. Consider $4$-subset $\{1,2,3,6\}$. It contains exactly two $3$-subsets from $C'$, namely $\{1,3,6\}$ and $\{2,3,6\}$ (see the definition of the code (M.2)).
The subset $\{1,2,3,6\}$ does not belong to $C$ because any vertex (i.e. $4$-subset) of the latter one contains $4-b_0=4$ $3$-subsets from $C'$. Therefore, $\{1,2,3,6\}$ is in $C_1$. However, for subset $\{1,2,3,5\}$ there is exactly one of its $3$-subset
in $C'$ and it is $\{2,3,5\}$. We see that $\{1,2,3,5\}$ is in $C_1$ and $b_1$ is not a constant, so $C$ in this case is not completely regular.

According to Corollary \ref{CoroJn3} there are several completely regular codes in $J(10,3)$ having the same intersection array $\{9;9\}$ as the code (M.2) but not isomorphic to it. We show that regardless the code being (M.2) or not, the corresponding completely regular code $C$ in $J(10,4)$ does not exist.
There are constants
$$b_0=0<b_1<b_2=4,$$
where $x\in C_i$ if and only if $x$ has $b_{i}$ vertices from $C'_{i-1}$ and $4-b_{i}$ vertices of $C'_{i}$, and $C'$ is a completely regular code in $J(10,3)$ with intersection array $\{9;9\}$.
Note that for $b_1<4$ and up to opposite code we have that  $b_1$ is $2$ or $3$. Let $\{1,2,3,4\}$ be in $C$. From the constant $b_0$ we have that all its 3-subsets: $\{1,2,3\}$, $\{2,3,4\}$, $\{1,2,4\}$, $\{1,3,4\}$ are in $C'$.
Consider the union of the $3$-subsets of $\{1,2,3,i\}$ for $i=4,\ldots,10$. It is not hard to see that these subsets form the closed neighborhood of $\{1,2,3\}$ in $J(10,3)$.
Depending on whether the vertex $\{1,2,3,i\}$  is in $C$ or in $C_1$, its $3$-subsets contain exactly $4$ or $b_1$ of vertices from $C'$, one of which is $\{1,2,3\}$.
Suppose exactly $r$ of vertices $\{1,2,3,i\}$, $i=4,\ldots,10$ are in $C$. Then there are exactly $3r+(b_1-1)(7-r)$ neighbors of $\{1,2,3\}$ in $J(10,3)$ that are in $C'$. On the other hand, from the intersection array $\{9;9\}$ of $C'$ and valency $21$ of the graph $J(10,3)$
we see that a vertex of $C'$ has exactly $12$ neighbors in $C'$. Therefore we obtain $3r+(b_1-1)(7-r)=12$ which does not hold for  $b_1=2$ or $3$ and nonegative integer $r$, a contradiction.

If $C'$ is (GP.3), i.e  such that $\{x:|x|=3, \mbox{ all elements of }$x$ \mbox{ are even or odd}\}$, then we see that
$C$ is (GP.3'),  i.e. $C=\{x:|x|=4, \mbox{ all elements of }$x$ \mbox{ are even or odd}\}$.

\eproof

In the end of this Section, we will focus on completely regular codes of strength $1$ in $J(n,4)$ without minimum eigenvalue in the spectrum. Although we have not obtained a full characterization of these codes, we derive some significant restrictions on their feasible intersection arrays.

\begin{lemma}\label{L3}
	Let $C$ be a completely regular code in $J(n,4)$ with covering radius $2$, the intersection matrix
	$$\left(%
	\begin{array}{ccc}
		4n-16-\beta_0 & \beta_0 &0\\
		\gamma_1& 4n-16-\gamma_1-\beta_1 & \beta_1\\
		0&\gamma_2 & 4n-16-\gamma_2\end{array}%
	\right)$$ and distance partition $\{C_0=C, C_1, C_2\}$. Then the following equalities hold: $$|C_0|=\gamma_1\gamma_2\frac{n(n - 1)(n - 2)(n - 3)}{24(\beta_0\beta_1+\beta_0\gamma_2 + \gamma_1\gamma_2)} ,\,|C_1|=\beta_0\gamma_2\frac{n(n - 1)(n - 2)(n - 3)}{24(\beta_0\beta_1+\beta_0\gamma_2 + \gamma_1\gamma_2)},\,|C_2|=\beta_0\beta_1\frac{n(n - 1)(n - 2)(n - 3)}{24(\beta_0\beta_1+\beta_0\gamma_2 + \gamma_1\gamma_2)}.$$
\end{lemma}
\bproof
Let us count the number of the edges between $C_0$ and $C_1$. On the one hand, it is equal to
$|C_0|\beta_0$, but on the other hand, it is equal to $|C_1|\gamma_1$. For edges between $C_1$ and $C_2$, we have that $|C_1|\beta_1=|C_2|\gamma_2$. Since $C$ has covering radius $2$, we also have that $|C_0|+|C_1|+|C_2|={n \choose 4}$. These three equations give us the required.
\eproof

For our purposes we need to use the following well-known local structure of the Johnson graphs.
\begin{prop}\label{prop:local}
	Let $x$ be a vertex of $J(n,w)$ and let $N(x)$ be the set of its neighbors. Then the graph induced on $N(x)$ by $J(n,w)$ is a $w\times (n-w)$-grid with edges between vertices in a common row or column. Each row of the grid together with $x$ forms a Delsarte clique in $J(n,w)$.
\end{prop}

\begin{theorem}\label{T6}
	Let $C$ be a completely regular code in $J(n,4)$ with covering radius $2$ and $Spec(C)=\{\theta_0(n,4),\theta_2(n,4), \theta_3(n,4) \}$. Then $C$ has the intersection matrix
	$$\left(%
	\begin{array}{ccc}
		4n-16-\beta_0 & \beta_0 &0\\
		-\frac{(\beta_0 - 2n + 2)(\beta_0 - 3n + 6)}{\beta_0 - \gamma_2}& \beta_0+\gamma_2-n-8 & \frac{(\gamma_2-2n+2)(\gamma_2-3n+6)}{\beta_0-\gamma_2}\\
		0&\gamma_2 & 4n-16-\gamma_2\end{array}%
	\right)$$
	for some distinct positive integers $\beta_0$ and $\gamma_2$.
	
\end{theorem}

\bproof
The completely regular code $C$ has the following $3 \times 3$ intersection matrix:

$$M=\left(%
\begin{array}{ccc}
	4n-16-\beta_0 & \beta_0 &0\\
	\gamma_1& 4n-16-\gamma_1-\beta_1 & \beta_1\\
	0&\gamma_2 & 4n-16-\gamma_2\end{array}%
\right).$$

Note that each row's sum equals the valency of $J(n,4)$, namely, $4n-16$.
By direct calculations we have the following equality for the characteristic polynomial of $M$:

$$ P(x) = det(xI-M) = x^3 + \psi_2x^2 + \psi_1x+ \psi_0,$$
where
\begin{equation}\label{Eq:0}
	\begin{aligned}
		\psi_2 & = \beta_0 + \beta_1 + \gamma_1 + \gamma_2 - 12n + 48,\\
		\psi_1 & = \beta_0\beta_1 + \beta_0\gamma_2 + \gamma_1\gamma_2 - 8\beta_0n - 8\beta_1n - 8\gamma_1n - 8\gamma_2n + 48n^2 + 32\beta_0 + 32\beta_1 + 32\gamma_1 + 32\gamma_2 - 384n + 768,\\
		\psi_0& = - 4\beta_0\beta_1n - 4\beta_0\gamma_2n - 4\gamma_1\gamma_2n + 16\beta_0n^2 + 16\beta_1n^2 + 16\gamma_1n^2 + 16\gamma_2n^2 - 64n^3 + 16\beta_0\beta_1 + 16\beta_0\gamma_2 + 16\gamma_1\gamma_2 \\
		& - 128\beta_0n - 128\beta_1n - 128\gamma_1n - 128\gamma_2n + 768n^2 + 256\beta_0 + 256\beta_1 + 256\gamma_1 + 256\gamma_2 - 3072n + 4096.
	\end{aligned}
\end{equation}

By the Theorem assumption we know that $Spec(M)=\{\theta_0(n,4),\theta_2(n,4), \theta_3(n,4) \}$. Therefore,
\begin{equation}\label{Eq:1}
	\psi_2= -\theta_0(n,4)-\theta_2(n,4)-\theta_3(n,4),
\end{equation}
\begin{equation}\label{Eq:2}
	\psi_1= \theta_0(n,4)\theta_2(n,4)+\theta_0(n,4)\theta_3(n,4)+\theta_2(n,4)\theta_3(n,4),
\end{equation}
\begin{equation}\label{Eq:3}
	\psi_0= -\theta_0(n,4)\theta_2(n,4)\theta_3(n,4),
\end{equation}
From (\ref{Eq:0}) and (\ref{Eq:1}) we immediately have that $\gamma_1 = -\beta_0 - \beta_1 - \gamma_2 + 5n - 8. $
Using this equality, (\ref{Eq:2}) can be simplified and written in the following way:
\begin{equation}\label{Eq:2.1}
	\beta_1(\beta_0-\gamma_2) =  (\gamma_2-2n+2)(\gamma_2-3n+6).
\end{equation}

If $\beta_0\neq \gamma_2$ then $\beta_1$ equals $\frac{(\gamma_2-2n+2)(\gamma_2-3n+6)}{\beta_0-\gamma_2}.$ Direct substitution of $\beta_1$ in (\ref{Eq:0}) and (\ref{Eq:3}) shows that the equality (\ref{Eq:3}) takes place for all values of $\beta_0$ and $\gamma_2$. Hence, the intersection matrix of the code is equal to
$$M_1=\left(%
\begin{array}{ccc}
	4n-16-\beta_0 & \beta_0 &0\\
	-\frac{(\beta_0 - 2n + 2)(\beta_0 - 3n + 6)}{\beta_0 - \gamma_2}& \beta_0+\gamma_2-n-8 & \frac{(\gamma_2-2n+2)(\gamma_2-3n+6)}{\beta_0-\gamma_2}\\
	0&\gamma_2 & 4n-16-\gamma_2\end{array}%
\right)$$
for some positive integers $\beta_0$ and $\gamma_2$.

Let us consider the case $\beta_0=\gamma_2$. Clearly, the equality (\ref{Eq:2.1}) implies that $\gamma_2=2n-2$ (the case $\gamma_2=3n-6$ contradicts to $\gamma_1>0$).
After substitutions $\beta_0=\gamma_2=2n-2$ in (\ref{Eq:0}) and (\ref{Eq:3}) equation (\ref{Eq:3}) can be written in the following form:
$$8(\beta_1 + \gamma_1 - n + 4)(n - 4)(n - 7)=0.$$
Consequently, for $\beta_0=\gamma_2$, the intersection matrix of the code is equal to
$$M_2=\left(%
\begin{array}{ccc}
	2n-14 & 2n-2 &0\\
	\gamma_1& 3n-12 & n-4-\gamma_1\\
	0&2n-2 & 2n-14\end{array}%
\right)$$
for some positive integer $\gamma_1$.

Our next goal is to prove that there are no completely regular codes with the intersection matrix $M_2$.
As usual, by $\{C_0=C, C_1, C_2\}$ we mean the distance partition of $C$. Lemma \ref{L3} guarantees us that
\begin{equation}\label{Eq:2.2}
	|C|=\frac{\gamma_1n(n-1)(n-3)}{72},\,|C_2|=\frac{(n-4-\gamma_1)n(n-1)(n-3)}{72}.
\end{equation}
By Proposition \ref{prop:local} we know that the graph induced by the set of the neighbors of an arbitrary vertex is $4 \times (n-4)$-grid with edges between vertices in a common row or column.

Without loss of generality we can consider
\begin{equation}\label{Eq:2.3}
 \gamma_1\leq \lfloor \frac{n-4}{2} \rfloor.
\end{equation}
 Take some $x\in C$ and let $N(x)$ be the grid corresponding to its neighbors. We know that $|N(x)\cap C_0|=2n-14$ and $|N(x)\cap C_1|=2n-2$. Since $2(n-4)>2n-14$, $N(x)$ may contain zero or exactly one row consisting of vertices from $C$ only.

\textbf{Case A.} Suppose that for every $y\in C$ we have that $N(y)$ contains exactly one row consisting of vertices from $C$ only. It means that we have a set of Delsarte cliques such that every code vertex belongs to exactly one of them. In other words, there exists a code $C'$ in $J(n,3)$ such that $C=\{x : |x|=4, \exists y\in C', y\subset x\}$, and code distance of $C'$ is at least $2$.

Since every vertex of $C'$ defines exactly $n-3$ vertices of $C$, we conclude that $C'=\frac{\gamma_1n(n-1)}{72}$. It is known that the size of a code in $J(n,3)$ with code distance $2$ is bounded above by $\frac{n(n-1)}{6}$, which follows from the Delsarte bound for independent sets. Therefore, we conclude that $\gamma_1\leq 12$. By our assumption for every $y\in C$ exactly one of $N(y)$'s rows consists of $n-4$ code vertices. Therefore, at least one of remaining three rows contains at least $\lceil \frac{n-10}{3} \rceil$ (and less than $n-4$) code vertices. Consequently, a vertex from $C_1$ from this row has at least $\lceil \frac{n-10}{3} \rceil+2$ neighbors from $C$, where one them is $y$ itself, and the other is taken from the row consisting of vertices from $C$ only. We have that
\begin{equation}\label{eq3}
	\lceil \frac{n-4}{3}\rceil \leq \gamma_1.
\end{equation}
In particular, $\gamma_1\leq 12$ implies $n\leq 40$. Since $N(y)$ contains $2n-14$ code vertices of $C$ and $n-4$ of them are located in one row, we conclude that $n\geq 10$.
As we noted before, $|C|=\frac{\gamma_1n(n-1)(n-3)}{72}$. By Proposition \ref{prop_strength} the code $C$ is a $1$-design. In other words, the number of the vertices of $C$ containing a fixed element from $\{1,2,\dots, n \}$ does not depend on the choice of the element. Clearly, this number equals $\frac{4}{n}|C|=\frac{\gamma_1(n-1)(n-3)}{18}$.
In above, we showed how $C$ is defined by $C'$. It is easy to see, that $C'$ also must be a $1$-design, so the number of the vertices of $C'$ containing a fixed element is equal to $\frac{3}{n}|C'|=\frac{\gamma_1(n-1)}{24}$.
Combining these arguments together, we have that
\begin{equation}\label{eq4}
	\frac{\gamma_1n(n-1)(n-3)}{72},\, \frac{\gamma_1(n-1)(n-3)}{18},\,\frac{\gamma_1(n-1)}{24} \in \mathbb{Z}.
\end{equation}
Consider the set $S=\{x\in C': 1\in x\}$. We know that $|S|=\frac{\gamma_1(n-1)}{24}$ and for any distinct $x$ and $y$ from $S$ we have that $|x\cap y|\leq 1$. Without loss of generality, $S=\big\{\{1,2i,2i+1\}, i=1,2, \dots \frac{\gamma_1(n-1)}{24}\big\}$.
Consider the set $S'=\{x\in C_2: 1\in x\}$. Since $C_2$ ia also a completely regular code of strength $1$, $|S'|=\frac{4}{n}|C_2|=\frac{(n-4-\gamma_1)(n-1)(n-3)}{18}$. Since the intersection size of any two vertices from $S$ and $S'$ is most $1$, for any $x\in S'$ we have that $x\cap \{2,3,\dots, \frac{\gamma_1(n-1)}{12}+1\}=\emptyset$. Consequently,
\begin{equation}\label{eq5}
	\frac{(n-4-\gamma_1)(n-1)(n-3)}{18}\leq {n-1-\frac{\gamma_1(n-1)}{12} \choose 3}.
\end{equation}
As we established above, $1\leq \gamma_1\leq 12$ and $10\leq n\leq 40$. Numerical verification of the fulfillment of conditions (\ref{eq4}) and (\ref{eq5}) for these values of $n$ and $\gamma_1$ gives only two possible cases:  $n=10$, $\gamma_1=8$ and $n=13$, $\gamma_1=6$. However, they both contradict to inequality (\ref{Eq:2.3}).

\textbf{Case B.} Suppose that for some $x\in C$, $N(x)$ does not have rows consisting of vertices from $C$ only (let us note, that this can be the case for other code vertices). Hence, there is a row with at least $\big\lceil \frac{2n-14}{4} \big\rceil$ code vertices. By our assumption this row also contains a vertex from $C_1$. Therefore, the vertex has at least $\big\lceil \frac{2n-14}{4} + 1 \big\rceil$ neighbors from $C_0$. As a result, we have that $$ \big\lceil \frac{n-5}{2} \big\rceil \leq \gamma_1 \leq \big\lfloor \frac{n-4}{2} \big\rfloor. $$

Analysis of these inequalities shows that it is possible in two cases:  $n$ is even, $\gamma_1=\frac{n-4}{2}$; and $n$ is odd, $\gamma_1=\frac{n-5}{2}$.
We need to put $2n-14$ vertices from $C$ in $4\times (n-4)$-grid such that any other vertex of the grid would have not more than $\frac{n-6}{2}$ neighbors from these vertices. Clearly, in the first case (i.e. when $n$ is even and $\gamma_1=\frac{n-4}{2}$) $N(x)$ must contain at least one row with $\frac{n-6}{2}$ code vertices. The remaining $\frac{3n-22}{2}$ code vertices must be located in other three rows in the same set of columns. In other words, $N(x)$ contains a $4\times(\frac{n-6}{2})$-grid with $2n-14$ and $2$ vertices from $C$ and $C_1$ respectively. There are only three nonequivalent ways for distribution of $C$ in a grid and for all of them one of these vertices from $C_1$ would have at least $\frac{n-4}{2}$ neighbors from $C$ inside this grid, a contradiction.
In the second case (i.e. when $n$ is odd and $\gamma_1=\frac{n-5}{2}$), the only possible configuration of code vertices in $N(x)$ is a $4\times (\frac{n-7}{2})$-grid. Without loss of generality, we have that $x=\{1,2,3,4\}$ and $$\big\{\{1,2,3,4\}\setminus\{i\}\cup \{j\}: i=1,2,3,4;\,j=5,6,\dots, \frac{n+1}{2} \big\}\subseteq C.$$
Consider a vertex $y=\{1,2,3,5\}$. Clearly, one of $N(y)$'s rows consists of $\frac{n-7}{2}$ code vertices and $\frac{n-1}{2}$ vertices from $C_1$. Consequently, code vertices in $N(y)$ have the same structure, so repeating the arguments we provided for $x$, for $y$ we have that
 $$\big\{\{1,2,3,5\}\setminus\{i\}\cup \{j\}: i=1,2,3,5;\,j=4,6,7,\dots, \frac{n+1}{2} \big\}\subseteq C.$$
Now we can repeat these arguments for all neighbors of $x$ from $C$, then for all their neighbors from $C$ and so on until then we have that the set $L$ of all $4$-subsets of a set $D=\{1,2,3,\dots, \frac{n+1}{2} \}$ is a subset of $C$ (conceptually, this part of the proof is close to the one from Theorem \ref{T5}). Let us note that all neighbors of vertices from $L$ outside $L$ must belong to $C_1$.

Suppose that there is another vertex $x'\in C\setminus L$ with the same type of $N(x')\cap C$ (without rows consisting of vertices from $C$ only). Consequently, there is one more set $L'$ of all $4$-subsets of a set $D'$ such that $|D'|=\frac{n+1}{2}$.
First of all, our previous arguments guarantee that $|D\cap D'|\leq 3$ (otherwise, $D$ must coincide with $D'$). If $|D\cap D'|=3$, a vertex $z=(D\cap D')\cup \{i\}$ for $i\in D\setminus D'$ has at least $2n-14+\frac{n-5}{2}$ neighbors from $C$. However, $z\in C$ and we get a contradiction. In case $|D\cap D'|=2$ we have that $z'=(D\cap D')\cup \{i\}\cup \{j\}$ for $i\in D\setminus D'$ and  $j\in D'\setminus D$ must be a vertex from $C_1$. Direct calculations show that $z'$ has at least $n-5$ neighbors from $L\cup L'$, and this is not possible. Since $D$ and $D'$ are both of size $\frac{n+1}{2}$, the last case here is $|D\cap D'|= 1$. In this case,
by direct calculations we have that $\{x \in C: D\cap D'\subseteq x\}$ consists of at least $2{\frac{n-1}{2} \choose 3}=\frac{(n-1)(n-3)(n-5)}{24}$ vertices. On the other hand, since $C$ is a $1$-design, the size of this set is bounded above by $\frac{n}{4}|C|=\frac{(n-1)(n-3)(n-5)}{36}$. We obtain a contradiction.

We proved that for all $y\in C\setminus L$ we have that $N(y)$ contains exactly one row consisting of vertices from $C$ only. Again we have a set of Delsarte cliques such that every code vertex (except the set $L$) belongs to exactly one of them. Therefore, there exists a code $C'$ in $J(n,3)$ such that $C\setminus L=\{x : |x|=4, \exists y\in C', y\subset x\}$, and code distance of $C'$ is at least $2$. By direct calculations and from (\ref{Eq:2.2}), we see that $|C'|=\frac{|C|-|L|}{n-3}=\frac{1}{n-3}(\frac{n(n-1)(n-3)(n-5)}{144}-\frac{(n+1)(n-1)(n-3)(n-5)}{384})=\frac{(5n-3)(n-1)(n-5)}{1152}$.
Since for any $x\in L$ we have that $N(x)\setminus L\subseteq C_1$, we conclude that for any $y\in C'$, $|y\cap D|=1$ or  $|y\cap D|=0$. Let us consider $S_1=\{y\in C':|y\cap D|=1\}$ and $S_0=\{y\in C':|y\cap D|=0\}$. Clearly,
\begin{equation}\label{eq6}
 |S_1|+|S_0|=|C'|=\frac{(5n-3)(n-1)(n-5)}{1152}.	
\end{equation}
Consider the following value:
$$Q=\sum_{t=\frac{n+3}{2}}^{n}{|\{y\in C: t\in y\}|}. $$
The fact that $C$ is a $1$-design implies that $Q=\frac{n-1}{2}\frac{4}{n}|C|=\frac{(n-1)^2(n-3)(n-5)}{72}$. The vertices from $L$ do not contribute to $Q$. Every clique defined by a vertex from $S_0$ adds $\frac{7n-25}{2}$ to this sum, and a clique defined by a vertex from $S_1$ adds $\frac{5n-17}{2}$.
Therefore, we have that $\frac{7n-25}{2}|S_0|+\frac{5n-17}{2}|S_1|=Q$. Together with (\ref{eq6}) it implies that $$|S_0|=\frac{(n - 1)(n - 5)(7n^2-28n+45)}{2304(n-4)}=\frac{1}{2304}\left(7n^3 - 42n^2 + 80n - 90 -\frac{135}{n-4}\right).$$
Since $|S_0|$ must be integer, $n-4$ divides $135$, so $n$ equals $9$, $13$, $19$, $31$, $49$ or $139$. Among them, $S_0$ is integral for $n=9$ and $n=13$, but $|C'|$ is non-integral for all of them. Finally, we proved that there are no completely regular codes with the intersection matrix $M_2$.
\eproof

Combining results from Theorems \ref{T5.5} and \ref{T6}, we formulate the main result for completely regular codes in $J(n,4)$.

\begin{theorem}\label{T7}
Let $C$ be a completely regular code of strength $1$ in $J(n,4)$. Then one of the following holds:
\begin{enumerate}
    \item $C$ has covering radius $1$, $n$ is $8$ and $C$ has the same intersection array as some of the codes from \cite{Vor} or \cite{AM2};
    \item $C$ has covering radius $2$, $Spec(C)=\{\theta_0(n,4),\theta_2(n,4), \theta_4(n,4) \}$ and $C$ is (GP.3') or (M.1) for $w=4$;
    \item $C$ has covering radius $2$ and $Spec(C)=\{\theta_0(n,4),\theta_2(n,4), \theta_3(n,4) \}$. Then $C$ has the intersection matrix
    $$\left(%
    \begin{array}{ccc}
    	4n-16-\beta_0 & \beta_0 &0\\
    	-\frac{(\beta_0 - 2n + 2)(\beta_0 - 3n + 6)}{\beta_0 - \gamma_2}& \beta_0+\gamma_2-n-8 & \frac{(\gamma_2-2n+2)(\gamma_2-3n+6)}{\beta_0-\gamma_2}\\
    	0&\gamma_2 & 4n-16-\gamma_2\end{array}%
    \right)$$
    for some distinct positive integers $\beta_0$ and $\gamma_2$.
\end{enumerate}

\end{theorem}

\bproof

The case $\theta_4(n,4) \in Spec(C)$ was considered in Theorem \ref{T5.5}. If $\theta_4(n,4) \not\in Spec(C)$ and $C$ has covering radius $2$, we use Theorem \ref{T6}. The last possible case is when $C$ has covering radius $1$ and $Spec(C)=\{\theta_0(n,4), \theta_2(n,4)\}$. In \cite{Vor} it was proved that such codes do not exist for $n\geq 9$. All possible intersection matrices for $n=8$ were listed in \cite{AM2} earlier.

\eproof

\section{Concluding remarks}
In this paper, we suggest an approach for tackling the existence problem of completely regular codes with minimum eigenvalue in geometric graphs. An approach showed promising results for the codes in Johnson graphs as we classified several series of  completely regular
codes in these graphs.

The method suggested in Section 3 can also work for completely regular codes in other geometric graphs.
From this perspective we view the Grassmann graphs as the most natural series for this technique. We suppose that some of the ideas from Section 4 can work for completely regular codes in Grassmann graphs with covering radius greater than $1$ with the minimum eigenvalue. In the case of binary field, the first open cases are completely regular codes in the Grassmann graphs $J_2(n,3)$ with covering radius $2$. From Section 3 we observe that these codes are in a natural correspondence with completely regular codes in $J_2(n,2)$. The latter codes for some spectra  were classified in \cite{FI} and one can use this restriction for classification of codes in $J_2(n,3)$.

In Section 4 we also described all the completely regular codes of strength $1$ in the Johnson graphs $J(n,4)$ with only one case for the eigenvalues $\theta_0(n,4), \theta_2(n,4), \theta_3(n,4)$ partially left open.
As one may see from the proof of Theorem \ref{T6} the intersection matrix $M_2$ depends on two parameters instead of one (it also holds for $M_1$). Although ideas that we used for proving the nonexistence of codes with the intersection matrix $M_1$ may be also helpful for the case with the matrix $M_2$ for some particular parameters, in our opinion, such codes might exist.

As an example, we present the following intersection matrices of possible completely regular codes with covering radius $2$ in $J(n,4)$ for small $n$ (all divisibility conditions known to us are satisfied here):

	 $$M=\left(%
	\begin{array}{ccc}
		14 & 30 & 0\\
		3& 31 & 10\\
		0&24 & 20\end{array}%
	\right) \text{ in } J(15,4),\,\,\,	
	M=\left(%
	\begin{array}{ccc}
		16 & 36 & 0\\
		3& 35 & 14\\
		0&24 & 28\end{array}%
	\right) \text{ in } J(17,4),\,\,\,	
	M=\left(%
	\begin{array}{ccc}
		20 & 36 & 0\\
		6& 42 & 8\\
		0&32 & 24\end{array}%
	\right) \text{ in } J(18,4).
	$$

{\it Acknowledgement.} The authors thank anonymous referee for thorough reading of the manuscript and substantial remarks that improved the presentation.

\end{document}